\documentclass[12pt]{article}

\usepackage{amssymb,amsmath,amsthm}
\usepackage{amsfonts}
\usepackage{enumerate}

\usepackage{authblk}

\usepackage{palatino}
\usepackage{fancyhdr}	
\usepackage{mdframed}

\usepackage{pgf,tikz}
\usepackage{subfig}
\usetikzlibrary{arrows}                      
\usetikzlibrary{decorations.markings}

\usepackage{geometry}
\newtheoremstyle{mytheoremstyle} % name
    {10pt}                    % Space above
    {10pt}                    % Space below
    {\normalfont}                   % Body font
    {}                           % Indent amount
    {\bfseries}                   % Theorem head font
    {.}                          % Punctuation after theorem head
    {0.3cm}                       % Space after theorem head
    {}  % Theorem head spec (can be left empty, meaning ‘normal’)
\theoremstyle{mytheoremstyle}

\theoremstyle{plain}% Theorem-like structures provided by amsthm.sty

\fancypagestyle{plain}{%
  \fancyhf{}%
}

\newtheorem{theorem}{Theorem}[section]

\newtheorem{lemma}[theorem]{Lemma}
\newtheorem{proposition}[theorem]{Proposition}

\newtheorem{remark}[theorem]{Remark}

\newtheorem{notation}[theorem]{Notation}

\begin{document}

\title{The simplicity index of tournaments}
\author{
	Abderrahim Boussa\"{\i}ri\thanks{Corresponding author: Abderrahim Boussaïri. Email: aboussairi@hotmail.com} , Soufiane Lakhlifi and Imane Talbaoui
}
\affil{Laboratoire Topologie, Alg\`ebre, G\'eom\'etrie et Math\'ematiques Discr\`etes, Facult\'e des Sciences A\"in Chock, Hassan II University of Casablanca, Maroc.}

\maketitle

\begin{abstract}
An $n$-tournament $T$ with vertex set $V$ is simple if there is no 
subset $M$ of $V$ such that $2\leq \left \vert M\right \vert \leq n-1$ 
and for every $x\in V\setminus M$, either $M\rightarrow x$ or $x
\rightarrow M$. The simplicity index of an $n$-tournament $T$ is the 
minimum number $s(T)$ of arcs whose reversal yields a non-simple 
tournament. M\"{u}ller and Pelant (1974) proved that $s(T)\leq\frac{n-1}{
2}$, and that equality holds if and only if $T$ is doubly regular. As 
doubly regular tournaments exist only if $n\equiv 3\pmod{4}$, $s(T)<\frac
{n-1}{2}$ for $n\not\equiv3\pmod{4}$. In this paper, we study the class 
of $n$-tournaments with maximal simplicity index for $n\not\equiv3\pmod{4
}$.
\end{abstract}



\section{Introduction}
  
	A \emph{tournament} $T$ consists of a finite set $V$ of \emph{vertices}
 together with a set $A$ of ordered pairs of distinct vertices, called 
\emph{arcs}, such that for all $x\not=y\in V$, $(x,y)\in A$ if and only 
if $(y,x)\notin A$. Such a tournament is denoted by $T=(V,A)$.  Given $x
\neq y\in V$,  we say that $x$ \emph{dominates} $y$ and we write $x
\rightarrow y$ when $(x,y)\in A$. Similarly, given two disjoint subsets $
X$ and $Y$ of $V$, we write $X\rightarrow Y$ if $x\rightarrow y$ holds 
for every $(x,y)\in X\times Y$.  Throughout this paper, we mean by an 
\emph{$n$-tournament} a tournament with $n$ vertices. 

     A tournament is \emph{regular} if there is an integer $k\geq 1$ such that each vertex dominates exactly $k$ vertices.
     It is \emph{doubly regular} if there is an integer $k\geq 1$ such that every unordered pair of vertices dominates exactly $k$ vertices.
    
    A tournament is \emph{transitive}, if for any vertices $x$, $y$ and $z$, $x\rightarrow y$ and $y\rightarrow z$ implies that $x\rightarrow z$. 
    A tournament $T=(V,A)$ is {\em reducible} if $V$ admits a bipartition $\{X,Y\}$ such that 
    $X\rightarrow Y$. 
    The notion of simple tournament was introduced by Fried and Lakser \cite{fried1971simple}, it was motivated by questions in algebra. It is closely related to modular decomposition \cite{gallai1967transitiv} which involves the notion of module. Recall that a \emph{module} of a tournament $T=(V,A)$ is a subset $M$ of $V$ such that for every $x\in V\setminus M$ either $M\rightarrow\left\{x\right\}$ or $\left\{x\right\}\rightarrow M$. For example, $\emptyset$, $\{x\}$, where $x\in V$, and $V$ are modules of $T$ called \emph{trivial} modules.  
    An $n$-tournament is \emph{simple} \cite{erdos1972some,muller1974strongly} (or prime \cite{cournier1992efficient} or primitive \cite{ehrenfeucht1990primitivity} or indecomposable \cite{ille1997indecomposable,schmerl1993critically}) if $n\geq 3$ and all its modules are trivial. The simple tournaments with at most $5$ vertices are shown in Figure \ref{fig1}. A tournament is \emph{decomposable} if it admits a non-trivial module.

\begin{figure}[htb]
\centering
\tikzset{fleche/.style args={#1:#2:#3:#4}{postaction = decorate, decoration={name=markings,mark=at position #1 with {\arrow[#2,scale=#3]{#4}}}}}
                                           
	\subfloat{ 
		\begin{tikzpicture}
			\def\r{2}
			\def\t{18}
			\coordinate (2) at (0+\t:\r) ;
			\coordinate (1) at (72+\t:\r) ;
			\coordinate (3) at (2*72+\t:\r) ;
																																		
			\draw[fill=black, color=black] (1) circle (1.4pt); 
			\draw[fill=black, color=black] (2) circle (1.4pt); 
			\draw[fill=black, color=black] (3) circle (1.4pt); 
																																		
			\node at (1) { } ; 
			\node at (2) { } ; 
			\node at (3) { } ;
			\draw[fleche=0.85:black:1.2:stealth,black,line width=1.2pt]  (1.center) to (2.center);
			\draw[fleche=0.85:black:1.2:stealth,black,line width=1.2pt]  (2.center) to (3.center);
			\draw[fleche=0.85:black:1.2:stealth,black,line width=1.2pt]  (3.center) to (1.center);
		\end{tikzpicture} 
	}
\qquad
	\subfloat{
		\begin{tikzpicture}
			\def\r{2}
			\def\t{18}
																			
			\coordinate (2) at (0+\t:\r) ;
			\coordinate (1) at (72+\t:\r) ;
			\coordinate (5) at (2*72+\t:\r) ;
			\coordinate (4) at (3*72+\t:\r) ;
			\coordinate (3) at (4*72+\t:\r) ;
																			
			\draw[fill=black, color=black] (1) circle (1.4pt); 
			\draw[fill=black, color=black] (2) circle (1.4pt); 
			\draw[fill=black, color=black] (3) circle (1.4pt); 
			\draw[fill=black, color=black] (4) circle (1.4pt);
			\draw[fill=black, color=black] (5) circle (1.4pt);

			\node   at (1) { } ; 
			\node  at (2) { } ; 
			\node at (3) { } ; 
			\node at (4) { } ;
			\node  at (5) { } ; 
																			
			\draw[fleche=0.85:black:1.2:stealth,black,line width=1.2pt]  (1.center) to (2.center);
			\draw[fleche=0.85:black:1.2:stealth,black,line width=1.2pt]  (2.center) to (5.center);
			\draw[fleche=0.85:black:1.2:stealth,black,line width=1.2pt]  (2.center) to (3.center);
			\draw[fleche=0.85:black:1.2:stealth,black,line width=1.2pt]  (3.center) to (4.center);
			\draw[fleche=0.85:black:1.2:stealth,black,line width=1.2pt]  (3.center) to (1.center);
			\draw[fleche=0.85:black:1.2:stealth,black,line width=1.2pt]  (5.center) to (3.center);
			\draw[fleche=0.85:black:1.2:stealth,black,line width=1.2pt]  (5.center) to (1.center);
			\draw[fleche=0.85:black:1.2:stealth,black,line width=1.2pt]  (4.center) to (2.center);
			\draw[fleche=0.85:black:1.2:stealth,black,line width=1.2pt]  (4.center) to (5.center);
			\draw[fleche=0.85:black:1.2:stealth,black,line width=1.2pt]  (1.center) to (4.center);
		\end{tikzpicture}
	}
\\
\vspace{0.3cm}

	\subfloat{
		\begin{tikzpicture}
			\def\r{2}
			\def\t{18}
																							
			\coordinate (2) at (0+\t:\r) ;
			\coordinate (1) at (72+\t:\r) ;
			\coordinate (5) at (2*72+\t:\r) ;
			\coordinate (4) at (3*72+\t:\r) ;
			\coordinate (3) at (4*72+\t:\r) ;

			\draw[fill=black, color=black] (1) circle (1.4pt); 
			\draw[fill=black, color=black] (2) circle (1.4pt); 
			\draw[fill=black, color=black] (3) circle (1.4pt); 
			\draw[fill=black, color=black] (4) circle (1.4pt);
			\draw[fill=black, color=black] (5) circle (1.4pt);
																							
			\node   at (1) { } ; 
			\node  at (2) { } ; 
			\node at (3) { } ; 
			\node at (4) { } ; 
			\node at (5) { } ; 
																							
			\draw[fleche=0.85:black:1.2:stealth,black,line width=1.2pt]  (1.center) to (2.center);
			\draw[fleche=0.85:black:1.2:stealth,black,line width=1.2pt]  (1.center) to (4.center);
			\draw[fleche=0.85:black:1.2:stealth,black,line width=1.2pt]  (2.center) to (3.center);
			\draw[fleche=0.85:black:1.2:stealth,black,line width=1.2pt]  (2.center) to (4.center);
			\draw[fleche=0.85:black:1.2:stealth,black,line width=1.2pt]  (2.center) to (5.center);
			\draw[fleche=0.85:black:1.2:stealth,black,line width=1.2pt]  (3.center) to (1.center);
			\draw[fleche=0.85:black:1.2:stealth,black,line width=1.2pt]  (3.center) to (4.center);
			\draw[fleche=0.85:black:1.2:stealth,black,line width=1.2pt]  (4.center) to (5.center);
			\draw[fleche=0.85:black:1.2:stealth,black,line width=1.2pt]  (5.center) to (1.center);
			\draw[fleche=0.85:black:1.2:stealth,black,line width=1.2pt]  (5.center) to (3.center);                               
		\end{tikzpicture}   
	}
\qquad
	\subfloat{
		\begin{tikzpicture}
			\def\r{2}
			\def\t{18}
			\coordinate (2) at (0+\t:\r) ;
			\coordinate (1) at (72+\t:\r) ;
			\coordinate (5) at (2*72+\t:\r) ;
			\coordinate (4) at (3*72+\t:\r) ;
			\coordinate (3) at (4*72+\t:\r) ;

			\draw[fill=black, color=black] (1) circle (1.4pt); 
			\draw[fill=black, color=black] (2) circle (1.4pt); 
			\draw[fill=black, color=black] (3) circle (1.4pt); 
			\draw[fill=black, color=black] (4) circle (1.4pt);
			\draw[fill=black, color=black] (5) circle (1.4pt);
																							
			\node   at (1) { } ; 
			\node  at (2) { } ; 
			\node at (3) { } ; 
			\node at (4) { } ; 
			\node  at (5) { } ; 
																							
			\draw[fleche=0.85:black:1.2:stealth,black,line width=1.2pt]  (1.center) to (2.center);
			\draw[fleche=0.85:black:1.2:stealth,black,line width=1.2pt]  (1.center) to (3.center);
			\draw[fleche=0.85:black:1.2:stealth,black,line width=1.2pt]  (1.center) to (4.center);
			\draw[fleche=0.85:black:1.2:stealth,black,line width=1.2pt]  (2.center) to (3.center);
			\draw[fleche=0.85:black:1.2:stealth,black,line width=1.2pt]  (2.center) to (5.center);
			\draw[fleche=0.85:black:1.2:stealth,black,line width=1.2pt]  (3.center) to (4.center);
			\draw[fleche=0.85:black:1.2:stealth,black,line width=1.2pt]  (3.center) to (5.center);
			\draw[fleche=0.85:black:1.2:stealth,black,line width=1.2pt]  (4.center) to (2.center);
			\draw[fleche=0.85:black:1.2:stealth,black,line width=1.2pt]  (5.center) to (1.center);
			\draw[fleche=0.85:black:1.2:stealth,black,line width=1.2pt]  (5.center) to (4.center);                                      
		\end{tikzpicture}
	}

\caption{The simple tournaments with at most $5$ vertices }
\label{fig1}   
\end{figure}

	Given an $n$-tournament $T$, the \emph{Slater index} $i(T)$ of $T$ is 
the minimum number of arcs that must be reversed to make $T$ transitive 
\cite{slater1961inconsistencies}. It is not difficult to see that $i(T)
\leq \frac{n(n-1)}{4}$. However, we do not know an exact determination 
of the upper bound of $i(T)$. Erd\H{o}s and Moon \cite{erdHos1965sets} 
proved that this bound is asymptotically equal to $\frac{n^{2}}{4}$.  
Recently, Satake \cite{satake2019constructive} proved that the Slater 
index of doubly regular $n$-tournaments is at least $\frac{n(n-1)}{4}-n^{
\frac{3}{2}}\log_2(2n)$.

Kirkland \cite{kirkland1993reversal} defined the {\em reversal index} $i_R(T)$
of a tournament $T$ as the minimum number of arcs whose reversal makes $T$ 
reducible. Clearly, $i_R(T)\leq i(T)$. Kirkland \cite{kirkland1993reversal} 
proved that $i_R(T)\leq \left \lfloor \frac{n-1}{2}\right \rfloor $ and 
characterized all the tournaments for which equality holds.

The indices above can be interpreted in terms of distance between 
tournaments. The \emph{distance} $d(T_1,T_2)$ between two tournaments $T_1$ 
and $T_2$ with the same vertex set is the number of pairs $\{x, y\}$ of 
vertices for which the arc between $x$ and $y$ has not the same direction in $
T_1$ and $T_2$. Let $\mathcal{F}$ be a family of tournaments with vertex set $
V$. The distance from a tournament $T$ to the family $\mathcal{F}$ is $d(T,
\mathcal{F})=\min\{ d(T,T^{\prime}):T^{\prime}\in \mathcal{F}\}$. If $\mathcal
{F}$ is the family of transitive tournaments on $V$, then $i(T)=d(T,\mathcal{F
})$. If $\mathcal{F}$ is the family of reducible tournaments on $V$, then $i_R
(T)=d(T,\mathcal{F})$.

By considering the family of decomposable tournaments, we obtain the simplicity 
index introduced by M\"{u}ller and Pelant \cite{muller1974strongly}. 
Precisely, consider an $n$-tournament $T$, where $n\geq 3$. The \emph{
simplicity index} $s(T)$ of $T$ (also called the \emph{arrow-simplicity} of $T
$ in \cite{muller1974strongly}) is the minimum number of arcs that must be 
reversed to make $T$ non-simple. For example, the tournaments shown in Figure 
\ref{fig1} have simplicity index $1$. Obviously, $s(T)\leq i_R(T)$ and $s(T)
\leq\frac{n-1}{2}$. M\"{u}ller and Pelant proved that $s(T)=\frac{n-1}{2}$ if 
and only if $T$ is doubly regular.

	A dual notion of the simplicity index is the decomposability index \cite{belkhechine2017decomposability}, 
which is 
obtained by considering the family of simple tournaments.

	In this paper, we provide an upper bound for $s(T)$, where $T$ is an $n
$-tournament for $n\not\equiv3\pmod{4}$. More precisely, we obtain the 
following result.

\begin{theorem}\label{MainTheorem}Given an $n$-tournament $T$, the following statements hold 
	\begin{enumerate}
    \item if $n=4k+2$, then $s(T)\leq2k$;
    \item if $n=4k+1$, then $s(T)\leq2k-1$;
    \item if $n=4k$, then $s(T)\leq2k-2$.
  \end{enumerate}
\end{theorem}
    
To prove that the bounds in this theorem are the best possible, we use 
the double regularity as follows.
    
\begin{theorem}
\label{TheoremTight}
  Let $l\in\{1,2,3\}$. Consider a doubly regular tournament $T$ of order 
$4k+3$, where $k\geq l$. The simplicity index of a tournament obtained 
from $T$ by removing $l$ vertices is $(2k+1)-l$.

\end{theorem}
     
    As shown by the next result, the opposite direction in Theorem~\ref{TheoremTight} holds when $l=1$. 
    \begin{theorem}
    \label{TheoremMain2}
  Given a tournament $T$ with $4k+2$ vertices, where $k\geq 1$, if $s(T)=2k$,
then $T$ is obtained from a doubly regular tournament by removing one vertex. 
    \end{theorem}
     
	The existence of doubly regular tournaments is equivalent to the 
existence of skew-Hadamard matrices \cite{reid1972doubly}. Wallis \cite{
wallis1971some} conjectured that $n\times n$ skew-Hadamard matrices 
exist if and only if $n=2$ or $n$ is divisible by $4$. Infinite 
families of skew-Hadamard matrices can be found in \cite{koukouvinos2008skew}.

	The most known examples of a doubly regular tournament are obtained 
from Paley construction.	For a prime power $q\equiv3\pmod{4}$, the \emph{Paley
tournament} of order $q$ is the tournament whose vertex set is the 
finite field $\mathbb{F}_q$, such that $x$ dominates $y$ if and only if $
x - y$ is a non-zero quadratic residue in $\mathbb{F}_q$.     

\section{Preliminaries}
    
  Let $T=(V,A)$ be an $n$-tournament and let $x\in V$. The \emph{
out-neighborhood} of $x$ is $N_{T}^{+}(x):=\left\{y\in V:x\rightarrow y
\right\}$, and the \emph{in-neighborhood} of $x$ is $N_{T}^{-}(x):=
\left\{y\in V:y\rightarrow x\right\}$. The \emph{out-degree} of $x$ (
resp. the \emph{in-degree} of $x$) is $\delta_{T}^{+}(x):=\left\vert N_{T
}^{+}(x)\right\vert$ (resp. $\delta_{T}^{-}(x):=\left\vert N_{T}^{-}(x)
\right\vert$). The \emph{out-degree} of $x$ is also called the \emph{
score} of $x$ in $T$. Recall that 
\begin{equation}
    \underset{z\in V}{{\displaystyle\sum}}\delta_{T}^{+}(z)=\underset{z\in V}{{\displaystyle \sum}}\delta_{T}^{-}(z)=\frac{n\left(n-1\right)}{2}.\label{SumDegrees}
\end{equation}
    A tournament is \emph{near-regular} if there exists an integer $k>0$ such that the out-degree of every vertex equals $k$ or $k-1$.
 
    \begin{remark}
    \label{RemarkRegularity} Let $T$ be an $n$-tournament. It follows from \eqref{SumDegrees} that
    \begin{enumerate}
    \item $T$ is regular if and only if $n$ is odd and every vertex has out-degree
    $\frac{(n-1)}{2}$;
    \item $T$ is near-regular if and only if $n$ is even and $T$ has $\frac{n}{2}$ vertices of out-degree $\frac{n}{2}$ and $\frac{n}{2}$ vertices of out-degree $\frac{(n-2)}{2}$.
    \end{enumerate}
    \end{remark}
    
  \begin{notation}\label{nota_PI}  
    Let $T=(V,A)$ be a near-regular tournament of order $4k+2$. We can partition $V$ into two $(2k+1)$-subsets, $V_{{\rm even}}:=\{z\in V,\delta_{T}^{+}(z)=2k\}$ and 
    $V_{{\rm odd}}:=\{z\in V,\delta_{T}^{+}(z)=2k+1\}$. 
    \end{notation}
     
     Let $x,y$ be distinct vertices of an $n$-tournament $T=(V,A)$. 
      The set
     $V\setminus\left\{x,y\right \}$ can be partitioned into four
     subsets:$\ N_{T}^{+}(x)\cap N_{T}^{+}(y)$, $N_{T}^{-}(x)\cap N_{T}^{-}(y)$, $N_{T}^{+}(x)\cap N_{T}^{-}(y)$ and $N_{T}^{-}(x)\cap
     N_{T}^{+}(y)$. The\emph{ out-degree} (resp. the \emph{in-degree}) of $(x,y)$ is
     $\delta_{T}^{+}(x,y):=\left \vert N_{T}^{+}(x)\cap N_{T}^{+}(y)\right \vert $ (resp. $\delta_{T}^{-}(x,y):=\left \vert N_{T}^{-}(x)\cap N_{T}^{-}(y)\right \vert $). The elements of $(N_{T}^{+}(x)\cap N_{T}^{-}(y))\cup(N_{T}^{-}(x)\cap N_{T}^{+}(y))$ are called \emph{separators }of $x,y$ and their
     number is denoted by $\sigma_{T}(x,y)$. 
    \begin{lemma}
    \label{od(xy)}Let $T$ be an $n$-tournament with vertex set $V$. For any $x\neq y\in V$, we have
    \begin{itemize}
    \item $\sigma_{T}(x,y)+\delta_{T}^{-}(x,y)+\delta_{T}^{+}(x,y)=n-2\label{SumSeparatorDegrees}$;
    \item $\delta_{T}^{-}(x,y)-\delta_{T}^{+}(x,y)=\delta_{T}^{-}(x)-\delta_{T}^{+}(y)$.
    \end{itemize}
    In particular, if $T$ is regular, then for any $x\neq y\in V$, $\delta_{T}^{-}(x,y)=\delta_{T}^{+}(x,y)$.
    \end{lemma}
    \begin{proof} 
    The first assertion is obvious. For the second assertion, we have $$\left \vert N_{T}^{-}(x)\right \vert =\left \vert N_{T}^{-}(x)\cap N_{T}^{-}(y)\right \vert
    +\left \vert N_{T}^{-}(x)\cap N_{T}^{+}(y)\right \vert +\left \vert N_{T}^{-}(x)\cap
    \left \{ y\right \} \right \vert $$
    and $$\left \vert N_{T}^{+}(y)\right \vert =\left \vert N_{T}^{+}(y)\cap N_{T}^{+}(x)\right \vert
    +\left \vert N_{T}^{+}(y)\cap N_{T}^{-}(x)\right \vert +\left \vert N_{T}^{+}(y)\cap
    \left \{ x\right \} \right \vert.$$
    Moreover, $y\in N_{T}^{-}(x)$ if and only if $x\in N_{T}^{+}(y)$. Then $\left \vert N_{T}^{-}(x)\cap \left \{  y\right \}  \right \vert =\left \vert N_{T}^{+}(y)\cap \left \{
    x\right \}  \right \vert $ and hence 
    \begin{equation*}
    \left \vert N_{T}^{-}(x)\cap N_{T}^{-}(y)\right \vert
    -\left \vert N_{T}^{+}(x)\cap N_{T}^{+}(y)\right \vert =\left \vert N_{T}^{-}(x)\right \vert -\left \vert N_{T}^{+}(y)\right \vert.\qedhere
    \end{equation*}
    \end{proof}
    Let $T=(V,A)$ be a tournament. For each vertex $z\in V$, we have
    \begin{equation*}
    \delta_{T}^{-}(z)\delta_{T}^{+}(z) = \left\vert \{\{x,y\}\in\tbinom{V}{2}:z\in (N_{T}^{-}(x)\cap N_{T}^{+}(y))\cup (N_{T}^{+}(x)\cap N_{T}^{-}(y))\}\right\vert.
    \end{equation*}
    By double-counting, we obtain 
    \begin{equation}
    \underset{z\in V}{\sum
    }\delta_{T}^{+}(z)\delta_{T}^{-}(z)=\underset{\left\{x,y\right\}\in \binom{V}{2}}{\sum}\sigma_{T}(x,y). \label{erdos1}%
    \end{equation}
      
    In the next proposition, we give some basic properties of doubly regular
    tournaments. For the proof, see \cite{muller1974strongly}.
    
    \begin{proposition}\label{DoublyRegular}
    Let $T=(V,A)$ be a doubly regular $n$-tournament. There exists $k\geq0$ 
    such that $n=4k+3$, $T$ is regular, and for all $x,y\in V$ such that $x\rightarrow y$, we have
    $$\left \vert N_{T}^{+}(x)\cap N_{T}^{+}(y)\right \vert =\left \vert N_{T}^{-}(x)\cap N_{T}^{-}(y)\right \vert =
    \left \vert N_{T}^{+}(x)\cap N_{T}^{-}(y)\right \vert =k$$ 
    $$\text{and}\hspace{3mm}\left \vert N_{T}^{-}(x)\cap N_{T}^{+}(y)\right \vert =k+1.$$
    \end{proposition}

\section{Proof of Theorem \ref{MainTheorem}}
    
    Let $T=(V, A)$ be a tournament. 
    Given a subset $B$ of $A$, we denote by ${\rm Inv}(T,B)$ the tournament obtained from $T$ by reversing all the arcs of $B$.
    We also use the following notation: $\delta_{T}^{+}$ $=\min \left \{  \delta_{T}^{+}(x):x\in V\right \}$, $\delta_{T}^{-}$ $=\min \left \{  \delta_{T}^{-}(x):x\in V\right \}$, $\delta_{T}=\min(\delta_{T}^{+},\delta_{T}^{-})$, and $\sigma_{T}=\min \{\sigma_{T}(x,y):x\neq y\in V\}$. 
    The next proposition provides an upper bound of the simplicity index of a tournament. 
    
    \begin{proposition}
    \label{Propositiondelta} For a tournament $T=(V,A)$ with at least $3$ vertices, we have
    $s(T)\leq \min(\delta_{T},\sigma_{T})$.
    \end{proposition}
    
    \begin{proof}
    Let $x\in V$. Clearly, the subset $V\setminus \left \{x\right \}  $ is a non-trivial module of 
    ${\rm Inv}(T,\left \{x\right \}\times N_{T}^{+}(x))$ and ${\rm Inv}(T,N_{T}^{-}(x)\times \left \{x\right \})$. It follows that $$s(T)\leq \min_{x\in V}(\delta_{T}^{+}(x),\delta_{T}^{-}(x))=\delta_{T}.$$
    
    Now, consider an unordered pair $\left \{x,y\right \}$ of vertices of $T$ and let
    $$B:=\left(\left \{x\right \}  \times\left((N_{T}^{+}(x)\cap N_{T}^{-}(y)\right)\cup\left(N_{T}^{+}(y)\cap N_{T}^{-}(x)\right)\times
    \left \{x\right \}\right).$$ Clearly, $\left \{x,y\right \}$ is a module of
    ${\rm Inv}(T,B)$. It follows that $$s(T)\leq \left \vert B\right \vert =\left \vert N_{T}^{+}(x)\cap N_{T}^{-}(y)\right \vert +\left \vert N_{T}^{+}(y)\cap N_{T}^{-}(x)\right \vert =\sigma_{T}(x,y).$$ Hence, $s(T)\leq \sigma_{T}$.
    \end{proof}
    In addition to the previous proposition, the proof of Theorem
    \ref{MainTheorem} requires the following lemma.
    \begin{lemma}
    \label{NmoinsUnSurDeux}Given an $n$-tournament $T=(V,A)$ with $n\geq$ $2$, we have 
    $$\text{$\delta_{T}\leq \left \lfloor \dfrac{n-1}{2}\right \rfloor $ and 
    $\ \sigma_{T}\leq \left \lfloor \dfrac{n-1}{2}\right \rfloor $.}$$
    \end{lemma}
    \begin{proof}
    For every $x\in V$, we have $\min\left(
    \delta_{T}^{+}(x),\delta_{T}^{-}(x)\right)  \leq \frac{n-1}{2}$. Thus, $$\delta_{T}\leq \left \lfloor \frac{n-1}{2}\right \rfloor .$$
    
    Now, to verify that $\sigma_{T}\leq \left \lfloor \dfrac{n-1}{2}
\right \rfloor $, observe that $$\sigma_{T}\leq\frac{1}{
\binom{\vert V \vert}{2}}\underset{\left \{  x,y\right \}  \in \binom{V}{
2}}{\sum}\sigma_{T}(x,y).$$ 
    It follows from \eqref{erdos1} that 
     \begin{align*}
    \sigma_{T} & \leq\frac{2}{n(n-1)}\underset{z\in V}{\sum}\delta_{T}^{+}%
    (z)\delta_{T}^{-}(z)\\
    & \leq\frac{2}{n(n-1)}\underset{z\in V}{\sum}\left(\frac{\delta_{T}%
    ^{+}(z)+\delta_{T}^{-}(z)}{2}\right)^{2}\\
    & \leq\frac{(n-1)}{2}.\qedhere
    \end{align*}
    \end{proof}
    
    \begin{proof}[Proof of Theorem \ref{MainTheorem}]
    For the first statement, suppose that $n=4k+2$. By Proposition \ref{Propositiondelta} and Lemma
    \ref{NmoinsUnSurDeux},\ we have $$s(T)\leq \delta_{T}\leq \left \lfloor \frac
    {n-1}{2}\right \rfloor =2k.$$
    
    For the second statement, suppose that $n=4k+1$. 
    By Proposition \ref{Propositiondelta}, $s(T)\leq\delta_{T}$. 
    If $T$ is not
    regular, then $\delta_{T}<\frac{n-1}{2}$ and hence $s(T)\leq2k-1$. Suppose that
    $T$ is regular and let $x\neq y\in V$. By Lemma \ref{od(xy)},
    $\sigma_{T}(x,y)=n-2-$ $\delta_{T}^{-}(x,y)-\delta_{T}^{+}(x,y)$ and $\delta_{T}^{-}(x,y)=\delta_{T}^{+}(x,y)$. Therefore, $\sigma_{T}(x,y)$ is odd, and hence $\sigma_{T}$ is odd as well. By Lemma \ref{NmoinsUnSurDeux}, $\sigma_{T}\leq \left \lfloor \frac{n-1}{2}\right \rfloor=2k$. Since $\sigma_{T}$ is odd, we obtain $\sigma_{T}\leq2k-1$. 
    It follows from Proposition \ref{Propositiondelta} that $s(T)\leq2k-1$.
    
    For the third statement, suppose that $n=4k$. 
    If $T$ is not near-regular, then $\delta_{T}<2k-1$, and hence $s(T)\leq2k-2$ 
    by Proposition~\ref{Propositiondelta}. Suppose that $T$ is near-regular.
    By Remark \ref{RemarkRegularity}, for every $z\in V$, $\delta_{T}^{+}%
    (z)\in \{2k,2k-1\}$. It follows from \eqref{erdos1} that%
    \begin{equation}
    \underset{\left \{  x,y\right \}  \in \binom{V}{2}}{\sum}\sigma_{T}%
    (x,y)=\underset{z\in V}{\sum}\delta_{T}^{+}(z)\delta_{T}^{-}(z)=8k^{2}%
    (2k-1).\label{SeparatorGlobal}%
    \end{equation}
    Thus, we obtain
    \begin{align*}
    \sigma_{T}&\leq\frac{1}{  \binom{\vert V \vert}{2} 
    }\underset{\left \{  x,y\right \}  \in \binom{V}{2}}{\sum}\sigma_{T}(x,y)\\
    &\leq\frac{2}{4k(4k-1)}8k^{2}(2k-1)\\
   &\leq(2k-1)+\frac{2k-1}{4k-1}\\
	 &\leq 2k-1.
    \end{align*}
    Since $s(T)\leq \sigma_{T}$ by Proposition \ref{Propositiondelta}, we obtain 
    $s(T)\leq \sigma_{T}\leq2k-1$. 
    Seeking a contradiction, suppose that $s(T)=2k-1$. We obtain $\sigma_{T}=2k-1$. 
    Let $x\in V_{{\rm even}}$ and $y\in V_{{\rm odd}}$ (see Notation~\ref{nota_PI}). It follows from Lemma \ref{od(xy)} that $\sigma_{T}(x,y)$ is
    even and hence $\sigma_{T}(x,y)\geq2k$. Thus, there are at least $(2k)^{2}$
    unordered pairs $\{x,y\}$ satisfying $\sigma_{T}(x,y)\geq2k$. For the other 
    $2\binom{2k}{2}$ unordered pairs, we have $\sigma_{T}(x,y)\geq \sigma_{T}=2k-1$. It follows that
    $$\underset{\left \{  x,y\right \}  \in \binom{V}{2}}{\sum}\sigma_{T}
    (x,y)\geq2\binom{2k}{2}(2k-1)+(2k)^{2}(2k)>8k^{2}(2k-1),$$ which contradicts
    (\ref{SeparatorGlobal}). Consequently, $s(T)\leq2k-2$.
    \end{proof}
    
    \section{Proof of Theorem \ref{TheoremTight}}
   
	To begin, recall that a \emph{graph} is defined by a vertex set $V$ 
and an edge set $E$. Two distinct vertices $x$ and $y$ of $G$ are \emph{
adjacent} if $\{x,y\}\in E$. For a vertex $x$ in $G$, the  set $N_G(x):=
\{y \in V: \{x,y\}\in E \}$ is called the \emph{neighborhood} of $x$ in $
G$. The \emph{degree} of $x$ is $\delta_G(x):=\left \vert N_G(x) \right 
\vert$. 
    
    Let $T=(V,A)$ be a tournament. To each subset $C$ of $V$, 
    we associate a graph in the following way. 
    Denote by $s_{C}(T)$ the minimum number of arcs that must be reversed to make
    $C$ a module of $T$. Clearly, 
    \begin{equation}\label{E_19_PI}
    s(T)=\min \left \{  s_{C}(T):2\leq \left \vert
    C\right \vert \leq n-1\right \}.
    \end{equation} 
    A graph $G=(V,E)$ is called a
    \emph{decomposability graph} for $C$ if $\left \vert E\right \vert =s_C(T)$ and $C$ is a module of the tournament $${\rm Inv}(T,\{(x,y)\in A:\{x,y\}\in E\})$$ 
    obtained from $T$ by reversing the arc between $x$ and $y$ for each edge 
    $\{x,y\}$ of $G$. 
        In the next lemma, we provide some of the properties of decomposability graphs.
    \begin{lemma}
    \label{LemmaGraphDec1} Let $T=(V,A)$ be a $n$-tournament and let $C$ be a
    subset of $V$ such that $2\leq \left \vert C\right \vert \leq n-1$. Given a decomposability graph $G=(V,E)$ for $C$, the following assertions hold
    \begin{itemize}
    \item $G$ is bipartite with bipartition $\left \{C,V\setminus C\right\}$;
    \item for each $x\in V\setminus C$, $N_{G}(x)=N_{T}^{+}(x)\cap C$ or $N_{G}(x)=N_{T}^{-}(x)\cap C$, and $\delta_{G}(x)=\min\left(\left\vert N_{T}^{-}(x)\cap C\right \vert ,\left \vert N_{T}^{+}(x)\cap C\right \vert \right)$.
    \end{itemize}
    \end{lemma}
    
    \begin{proof}
   The first assertion follows from the minimality of $\left \vert E\right \vert =s_C(T)$. 
        For the second assertion, consider $x\in V\setminus C$. 
        Since $C$ is a module of the tournament 
        ${\rm Inv}(T,\{(x,y)\in A:\{x,y\}\in E\})$, we have 
        $N_{G}(x)=N_{T}^{+}(x)\cap C$ or $N_{G}(x)=N_{T}^{-}(x)\cap C$. 
        Furthermore, it follows from the minimality of $\left \vert E\right \vert =s_C(T)$ that 
        $\delta_{G}(x)=\min\left(\left \vert N_{T}^{-}(x)\cap C\right \vert ,\left \vert N_{T}^{+}(x)\cap C\right \vert \right)$. 
    \end{proof}
    
    The next proposition is useful to prove Theorems \ref{TheoremTight} and \ref{TheoremMain2}. 
    
    \begin{proposition}
    \label{prop graphe decomp}
    Let $T=(V,A)$ be an $n$-tournament and let $C$ be a
    subset of $V$ such that $2\leq \left \vert C\right \vert \leq n-1$. Given a decomposability graph $G=(V,E)$ for $C$, the following statements hold
    \begin{itemize}
    \item if $n-\delta_{T}\leq \left \vert C\right \vert $, then $s_{C}(T)$
    $\geq \delta_{T}$;
    \item if $\left \vert C\right \vert \leq \sigma_{T}$, then $s_{C}(T)\geq \sigma_{T}$.
    \end{itemize}
    \end{proposition}
    \begin{proof}
    Before showing the first assertion, we establish 
    \begin{equation}\label{E_6a_PI}
    \left \vert E\right \vert  \geq(n-\left \vert C\right \vert )(\left \vert C\right \vert -(n-1-\delta_{T})).
    \end{equation}
    Let $x\in V\smallsetminus C$. By the second assertion of Lemma
    \ref{LemmaGraphDec1} 
    \begin{align*}
   \delta_{G}(x)&=\min ( \left \vert N_{T}^{-}(x)\cap C\right \vert ,\left \vert
    N_{T}^{+}(x)\cap C\right \vert )\\
    &=\left \vert C\right \vert -\max ( \left \vert N_{T}^{-}(x)\cap C\right \vert
    ,\left \vert N_{T}^{+}(x)\cap C\right \vert ).
    \end{align*}
    Therefore, we obtain  
     \begin{align}\label{E_7_PI}
    \delta_{G}(x)&\geq\left \vert C\right \vert -\max (\left \vert N_{T}^{-}(x)\right \vert,\left \vert N_{T}^{+}(x)\right \vert )\nonumber\\
    &
     \geq(\left \vert C\right \vert -(n-1-\delta_{T})).
    \end{align}
    Since $G$ is bipartite with bipartition $\left \{C,V\setminus C\right\}$, we have 
    \begin{equation*}
    \left \vert E\right \vert =
    {\displaystyle \sum \limits_{x\in V\smallsetminus C}}
    \delta_{G}(x).
    \end{equation*} 
    It follows from \eqref{E_7_PI} that 
    \begin{align*}
    \left \vert E\right \vert  & \geq \left \vert V\smallsetminus C\right \vert
    (\left \vert C\right \vert -(n-1-\delta_{T}))\\
    & \geq (n-\left \vert C\right \vert )(\left \vert C\right \vert -(n-1-\delta_{T})).
    \end{align*}
    Thus, \eqref{E_6a_PI} holds. 
    Moreover, we have $$(n-\left \vert C\right \vert )(\left \vert C\right \vert -(n-1-\delta
    _{T}))-\delta_{T}=(n-1-\left \vert C\right \vert )(\left \vert C\right \vert
    -(n-\delta_{T})).$$
    
    Now, to prove the first assertion, suppose that 
    $n-\delta_{T}\leq \left \vert C\right \vert $. 
    We obtain 
    $(n-1-\left \vert C\right \vert )(\left \vert C\right \vert
    -(n-\delta_{T})\geq0$, and hence $(n-\left \vert C\right \vert )(\left \vert C\right \vert -(n-1-\delta
    _{T}))\geq
    \delta_{T}$. It follows that
    $s_{C}(T)=\left \vert E\right \vert \geq
    \delta_{T}$. 
    
    Before showing the second assertion, we establish 
    \begin{equation}\label{E_6b_PI}
   \left \vert E\right \vert\geq \frac{\left\vert C\right\vert}{2}(2-\left\vert C\right\vert +\sigma_{T}).
    \end{equation}
    Consider two vertices $x\neq y\in C$. 
    For convenience, denote by $\mathcal{S}_T(x,y)$ the set of separators of 
    $\{x,y\}$. 
    Clearly, we have $(\mathcal{S}_T(x,y)\setminus C)\subseteq N_G(x)\cup N_G(y)$. 
    It follows that $$\delta_{G}(x)+\delta_{G}(y)\geq |\mathcal{S}_T(x,y)\setminus C|\geq
    \sigma_{T}(x,y)-(\left \vert C\right \vert -2).$$
   Consequently, we obtain 
    \begin{equation}
    \delta_{G}(x)+\delta_{G}(y)\geq \sigma_{T}-\left \vert C\right \vert +2. \label{eq6}
    \end{equation}
    Furthermore, observe that
     $$\sum_{\{x,y\}\in \binom{C}{2}}(\delta_{G}(x)+\delta_{G}(y))=(\left\vert
    C\right\vert -1)\sum_{x\in C}\delta_{G}(x).$$
    It follows from \eqref{eq6} that 
    $$(\left\vert C\right\vert -1)\sum_{x\in C}\delta_{G}(x)\geq 
    \dbinom{\left\vert C\right\vert }{2}(2-\left\vert C\right\vert +\sigma_{T}).$$
    Therefore, we have 
    $$\sum_{x\in C}\delta_{G}(x)\geq \frac{\left\vert C\right\vert}{2}(2-\left\vert C\right\vert +\sigma_{T}).$$
    Since $G$ is bipartite with bipartition $\left \{C,V\setminus C\right\}$, we have 
    \begin{equation*}
    \left \vert E\right \vert =\sum_{x\in C}\delta_{G}(x). 
    \end{equation*}
    We obtain 
    $$\left \vert E\right \vert\geq \frac{\left\vert C\right\vert}{2}(2-\left\vert C\right\vert +\sigma_{T}),$$ so \eqref{E_6b_PI} holds. 
    
    Finally, to prove the second assertion, suppose that 
    $\left\vert C\right\vert \leq \sigma_{T}$. 
    We obtain 
    $\frac{\left\vert C\right\vert}{2}(2-\left\vert C\right\vert +\sigma_{T})\geq \sigma_{T}$. 
    Since $s_{C}(T)=\left\vert E\right\vert$, it follows from \eqref{E_6b_PI} that $s_{C}(T)\geq \sigma_{T}$. 
    \end{proof}
    
    \begin{proof}[Proof of Theorem \ref{TheoremTight}] 
    Let $l\in\{1,2,3\}$. 
    Consider  a tournament $R$ from $T$ by removing $l$ vertices $v_1,\ldots,v_l$. 
    Set $V':= V\smallsetminus \{v_1,\ldots,v_l\}$. 
    It follows from Theorem \ref{MainTheorem} that $s(R)\leq (2k+1)-l$. 
    It remains to show that $s(R)\geq (2k+1)-l$. 
    By \eqref{E_19_PI}, it suffices to verify that $s_{C}(R)\geq (2k+1)-l$ for every 
    subset $C$ of $V'$ such that $2\leq \left \vert C\right \vert \leq (4k+2)-l$. 
    Let $C\subseteq V'$ such that $2\leq \left \vert C\right \vert \leq (4k+2)-l$. 
    We distinguish the following three cases. 
    \begin{itemize}
    \item Suppose that $2\leq \left \vert C\right \vert \leq (2k+1)-l$. 
    Since $T$ is doubly regular, it follows from Proposition \ref{DoublyRegular} that 
    $\sigma_{T}=2k+1$. 
    Therefore, $\sigma_{R}\geq (2k+1)-l$. 
    Since $2\leq \left \vert C\right \vert \leq (2k+1)-l$, 
    $\sigma_{R}\geq \left \vert C\right \vert$. 
    It follows from Proposition \ref{prop graphe decomp}
    that $s_{C}(R)\geq \sigma_{R}$, and hence $s_{C}(R)\geq (2k+1)-l$. 
    \item Suppose that $2k+2\leq \left \vert C\right \vert \leq (4k+2)-l$. 
    Since $T$ is doubly regular, it follows from Proposition \ref{DoublyRegular} that $T$ is regular. 
    Thus, $\delta_{T}=2k+1$. It follows that $\delta_{R}\geq (2k+1)-l$. 
    Since $2k+2\leq \left \vert C\right \vert \leq (4k+2)-l$, we obtain 
    $\left \vert C\right \vert + \delta_{R}\geq (4k+3)-l$. 
    It follows from Proposition \ref{prop graphe decomp}
    that $s_{C}(R)\geq \delta_{R}$, and hence $s_{C}(R)\geq (2k+1)-l$. 
    \item Suppose that $(2k+2)-l\leq \left \vert C\right \vert \leq 2k+1$. 
    Let $G=(E',V')$ be a decomposability graph for $C$. 
    We verify that 
    \begin{equation}\label{E_11_PI}
    |\{x\in V'\smallsetminus C:\delta_{G}(x)\neq 0\}|\geq |V'\smallsetminus C|-1.
    \end{equation}
    Otherwise, there exist $x\neq y\in V'\smallsetminus C$ such that 
    $\delta_{G}(x)=\delta_{G}(y)=0$. 
    It follows from the second assertion of Lemma \ref{LemmaGraphDec1} applied to $R$ that $C$ is contained in  one of the following intersections: $(N_{R}^{-}(x)\cap N_{R}^{+}(y))$, $(N_{R}^{-}(x)\cap N_{R}^{-}(y))$, $(N_{R}^{+}(x)\cap N_{R}^{+}(y))$, or 
    $(N_{R}^{+}(x)\cap N_{R}^{-}(y))$. 
    Thus, $C$ is contained in $(N_{T}^{-}(x)\cap N_{T}^{+}(y))$, $(N_{T}^{-}(x)\cap N_{T}^{-}(y))$, $(N_{T}^{+}(x)\cap N_{T}^{+}(y))$, or 
    $(N_{T}^{+}(x)\cap N_{T}^{-}(y))$. 
    It follows from Proposition \ref{DoublyRegular} that $\left \vert C\right \vert \leq k+1$, 
    which contradicts 
    $\left \vert C\right \vert \geq (2k+2)-l$ because $k\geq l$. 
    Consequently, \eqref{E_11_PI} holds. 
    Since $G$ is bipartite with bipartition $\left \{C,V'\smallsetminus C\right\}$, we have 
    $\left \vert E'\right \vert =\sum_{x\in V'\smallsetminus C}\delta_{G}(x)$. 
    Since $\left \vert E'\right \vert =s_{C}(R)$, we obtain 
        \begin{align*}
    s_{C}(R)&=\sum_{x\in V'\smallsetminus C}\delta_{G}(x)\\
    &\geq \left \vert V' \smallsetminus C \right \vert -1\hspace{5mm}\text{(by \eqref{E_11_PI})}\\
    &\geq (2k+1)-l\hspace{5mm}\text{(because $\left \vert C\right \vert \leq 2k+1$)}.\qedhere
    \end{align*}
    \end{itemize}
    \end{proof}

\section{Proof of Theorem \ref{TheoremMain2}}
    If a tournament $T$ is obtained from a doubly regular $(4k+3)$-tournament by deleting one vertex, then $T$ is near-regular and it follows from Proposition \ref{DoublyRegular} that
    
    \begin{itemize}
    \item[(C1)] if $x,y\in V_{{\rm even}}$ (see Notation~\ref{nota_PI}) or $x,y\in V_{{\rm odd}}$, then $\sigma_{T}(x,y)=2k+1$.
    \item[(C2)] if $x\in V_{{\rm even}}$ and $y\in V_{{\rm odd}}$, then $\sigma_{T}(x,y)=2k$.
    \end{itemize}
    Conversely, we have the following proposition. 
    \begin{proposition}
    \label{PropoLakhlifi}
    Let $T=(V,A)$ be a near-regular tournament of order $4k+2$. 
    If $T$ satisfies (C1) and (C2), then the tournament $U$ obtained from $T$ by adding a vertex $\omega$ which dominates $V_{{\rm odd}}$ and is dominated by $V_{{\rm even}}$ is doubly regular.
    \end{proposition}
     The proof of this proposition uses the following lemma.
    \begin{lemma}
    \label{LemLakhlifi}Under the notation and conditions of Proposition \ref{PropoLakhlifi}, for every $x,y\in V$ such that $x\rightarrow y$, we have
    \begin{itemize}
    \item if $x,y\in V_{{\rm odd}}$, then 
    $$\text{$\left \vert N_{T}^{-}(x)\cap N_{T}^{+}(y)\right \vert =\allowbreak
    k+1$ and $\left\vert N_{T}^{+}(x)\cap N_{T}^{-}(y)\right\vert=\allowbreak k$;}$$
    \item if $x,y\in V_{{\rm even}}$, then $$\text{$\left\vert N_{T}^{-}(x)\cap N_{T}^{+}(y)\right\vert =\allowbreak k+1$ and $\left \vert N_{T}^{+}(x)\cap N_{T}^{-}(y)\right \vert =\allowbreak k$;}$$
    \item if $x\in V_{{\rm odd}}$ and $y\in V_{{\rm even}}$, then $$\text{$\left\vert N_{T}^{-}(x)\cap N_{T}^{+}(y)\right \vert=\allowbreak k$ and $\left\vert N_{T}^{+}(x)\cap N_{T}^{-}(y)\right\vert=\allowbreak k$;}$$
    \item if $x\in V_{{\rm even}}$ and $y\in V_{{\rm odd}}$, then $$\text{$\left\vert N_{T}^{-}(x)\cap N_{T}^{+}(y)\right\vert=\allowbreak k+1$ and $\left\vert N_{T}^{+}(x)\cap N_{T}^{-}(y)\right\vert
    =\allowbreak k-1$.}$$
    \end{itemize}
    \end{lemma}
    \begin{proof}
    We have
    \begin{equation}\label{E_PI_0}
    \begin{cases}
    \left \vert N_{T}^{-}(x)\cap N_{T}^{-}(y)\right \vert +\left \vert N_{T}^{-}(x)\cap N_{T}^{+}(y)\right \vert
    =\left \vert N_{T}^{-}(x)\right \vert\\ 
    \text{and}\\ 
    \left \vert N_{T}^{+}(x)\cap N_{T}^{+}(y)\right \vert +\left \vert N_{T}^{-}(x)\cap N_{T}^{+}(y)\right \vert
    =\left \vert N_{T}^{+}(y)\right \vert.
    \end{cases}
    \end{equation}
    By using Lemma \ref{od(xy)}, we obtain 
    \begin{equation}\label{E_PI_1}
    \left \vert N_{T}^{-}(x)\cap N_{T}^{+}(y)\right \vert =\frac{1}{2}\left(\left \vert N_{T}^{-}(x)\right \vert+\left \vert N_{T}^{+}(y)\right \vert -4k+\sigma_{T}(x,y)\right).
     \end{equation}
     
    Using Assertions (C1) and (C2), we obtain the desired values of $ \vert N_{T}^{-}(x)\cap N_{T}^{+}(y)\vert$. Then, $ \vert N_{T}^{+}(x)\cap N_{T}^{-}(y) \vert$ follows immediately because
    $ \vert N_{T}^{+}(x)\cap N_{T}^{-}(y) \vert = \sigma(x,y)-\vert N_{T}^{-}(x)\cap N_{T}^{+}(y) \vert$.
    \end{proof}
    
    \begin{proof}[Proof of Proposition \ref{PropoLakhlifi}]
     Clearly, $U$ is regular. Furthermore, by Lemma
    \ref{od(xy)}, $$\delta_{U}^{+}(x,y)=\frac{4k-\sigma_{U}(x,y)+1}{2}$$ for distinct $x,y\in V\cup\{ \omega \}$. Therefore, $U$ is doubly regular if and only if $\sigma_{U}(x,y)=2k+1$ for every $x,y\in V\cup \{ \omega \}$. This equality follows directly from (C1) and (C2) when $x,y\in V$. Hence, it remains to prove that 
    \begin{equation}\label{E_5_P1}
    \text{$\sigma_{U}
    (\omega,z)=2k+1$ for every $z\in V$.} 
    \end{equation}
    Consider $z\in V$. It is not difficult to see that
    \begin{equation*}
    \sigma
        _{U}(\omega,z)=\left \vert N_{T}^{+}(z)\cap V_{{\rm even}}\right \vert +\left \vert N_{T}^{-}(z)\cap
        V_{{\rm odd}}\right \vert\mbox{ (see Notation~\ref{nota_PI})}.
    \end{equation*}
    
     Let $A_{{\rm odd}}:=(N_{T}^{+}(z)\cap V_{{\rm odd}})$,
    $A_{{\rm even}}:=( N_{T}^{+}(z)\cap V_{{\rm even}})$, $B_{{\rm odd}}:=(N_{T}^{-}(z)\cap
    V_{{\rm odd}})$, and $B_{{\rm even}}:=(N_{T}^{-}(z)\cap V_{{\rm even}})$. 
    We determine 
    $\left \vert A_{{\rm odd}}\right \vert $, $\left \vert A_{{\rm even}}\right \vert $, $\left \vert B_{{\rm odd}}\right \vert $, and $\left \vert B_{{\rm even}}\right \vert $ as follows. 
    
    To begin, suppose that
    $z\in V_{{\rm odd}}$. By counting the number of arcs from $N_{T}^{+}(z)$ to $N_{T}^{-}(z)$ in two ways, we get 
    \begin{align*}
    &\sum_{t\in A_{{\rm odd}}}
    \left \vert N_{T}^{-}(z)\cap N_{T}^{+}(t)\right \vert +
    \sum_{t\in A_{{\rm even}}}
    \left \vert N_{T}^{-}(z)\cap N_{T}^{+}(t)\right \vert \\
    =&\sum_{t\in B_{{\rm odd}}}
    \left \vert N_{T}^{-}(t)\cap N_{T}^{+}(z)\right \vert +
    \sum_{t\in B_{{\rm even}}}
    \left \vert N_{T}^{-}(t)\cap N_{T}^{+}(z)\right \vert .
    \end{align*}
    It follows from Lemma \ref{LemLakhlifi} that
    $$(k+1)\left \vert A_{{\rm odd}}\right \vert +k\left \vert A_{{\rm even}}\right \vert
    =(k+1)(\left \vert B_{{\rm odd}}\right \vert +\left \vert B_{{\rm even}}\right \vert).$$
    Since $z\in V_{{\rm odd}}$, we have $\left \vert A_{{\rm odd}}\right \vert +\left \vert
    A_{{\rm even}}\right \vert =2k+1$, $\left \vert B_{{\rm odd}}\right \vert +\left \vert
    B_{{\rm even}}\right \vert =2k$, $\left \vert A_{{\rm odd}}\right \vert +\left \vert B_{{\rm odd}}
    \right \vert =2k$, and $\left \vert A_{{\rm even}}\right \vert +\left \vert B_{{\rm even}}\right \vert
    =2k+1$.
    It follows that $\left \vert A_{{\rm odd}}\right \vert =\allowbreak k$, $\left \vert
    B_{{\rm odd}}\right \vert =k$, $\left \vert B_{{\rm even}}\right \vert =k$, and $\left \vert
    A_{{\rm even}}\right \vert =k+1$.
    
    Similarly, if $z\in V_{{\rm even}}$, then $\left \vert
    A_{{\rm odd}}\right \vert =\allowbreak k$, $\left \vert B_{{\rm odd}}\right \vert =k+1$,
    $\left \vert B_{{\rm even}}\right \vert =k$, and $\left \vert A_{{\rm even}}\right \vert =k$.
    
    Consequently, \eqref{E_5_P1} holds whatever the parity of $\delta^{+}_{T}(z)$.
    \end{proof}

    \begin{proof}[Proof of Theorem \ref{TheoremMain2}]
    Given $k\geq 1$, consider a tournament $T$, with $4k+2$ vertices, such that $s(T)=2k$. 
    By Proposition \ref{Propositiondelta}, $\delta_{T}\geq 2k$. 
   Thus, 
    $T$ is near-regular. 
    We conclude by applying Proposition \ref{PropoLakhlifi}. 
    Therefore, it suffices to verify that (C1) and (C2) are satisfied. 
    
    By Proposition \ref{Propositiondelta}, $\sigma_{T}(x,y)\geq2k$ for distinct 
    $x,y\in V$. 
    Moreover, it follows from Lemma \ref{od(xy)} that if $x,y\in V_{{\rm even}}$ or
    $x,y\in V_{{\rm odd}}$ (see Notation~\ref{nota_PI}), then $\sigma_{T}(x,y)$ is odd and hence $\sigma_{T}(x,y)\geq2k+1$. 
    
    Lastly, seeking a contradiction, suppose that (C1) or (C2) are not satisfied. 
   One of the following situations occurs 
    \begin{itemize}
    \item there are distinct $x,y\in V_{{\rm even}}$ such that $\sigma_{T}(x,y)>2k+1$,
    \item there are distinct $x,y\in V_{{\rm odd}}$ such that $\sigma_{T}(x,y)>2k+1$,
    \item there are $x\in V_{{\rm even}}$ and $y\in V_{{\rm odd}}$ such that $\sigma
    _{T}(x,y)>2k$.
    \end{itemize}
    We obtain 
    \begin{align*}
    \sum\limits_{\{x,y\}\in\binom{V}{2}}\sigma_{T}(x,y)&> (2k+1)\binom{\left \vert V_{{\rm even}}\right \vert }2+(2k+1)\binom{\left \vert V_{{\rm odd}}\right \vert }{2}+2k\left \vert V_{{\rm even}}\right \vert \left \vert V_{{\rm odd}}\right \vert\\
    &=4k(2k+1)^{2},
\end{align*}
    which contradicts (\ref{erdos1}). 
    Consequently, (C1) and (C2) are satisfied.
    \end{proof}

\section{Concluding remarks}

    \textbf{1.} An $n$-tournament with $n=4k+1$ is called \emph{near-homogeneous} 
    \cite{TABIB198077} if every unordered pair of its vertices belongs to $k$ or $(k+1)$ 3-cycles. 
    The existence of near-homogeneous tournaments is discussed in \cite{TABIB198077}, \cite{astie1992near}, and \cite{moukouelle1998construction}. 
    For $n\equiv1\pmod{4}$ or $n\equiv0\pmod{4}$, the $n$-tournaments given in 
    Theorem~\ref{TheoremTight} are not the only ones with a maximal simplicity index. Indeed, 
    let $T$ be a near-homogeneous tournament $T$ with $4k+1$ vertices. 
    By adapting the proof of Theorem \ref{TheoremTight}, we can verify that $s(T)=2k-1$. Moreover, by removing one vertex from $T$, we obtain a $(4k)$-tournament whose simplicity index is $2k-2$. 
    Consequently, an analogue of Theorem~\ref{TheoremMain2} does not exist when 
    $l=2$ or $3$. 
    
	\textbf{2.} The score vector of a $n$-tournament $T$ is the ordered sequence
of the scores of $T$ listed in a non-decreasing order. Kirkland
\cite{kirkland1993reversal} proved that the reversal index of an $n$-tournament $T$ is equal to
$$\min\left\{\sum_{i=1}^{j}s_{i}-\binom{j}{2}:{1\leq j\leq n}\right\},$$
where $(s_{1},s_{2},\ldots,s_{n})$ is the score vector of $T$.

    An equivalent form of this result was obtained earlier by Li and 
Huang \cite{li1989score}. As a consequence, two tournaments with the 
same score vector have the same reversal index. This fact is not true 
for the simplicity index. Indeed, for an odd number $n$, consider the
$n$-tournament $R_n$ whose vertex set is the additive group $\mathbb{Z}_n=
\{0,1,\ldots,n-1\}$ of integers modulo $n$, such that $i$ dominates $j$ if 
and only if $i-j \in\{1,\ldots,(n-1)/2\}$. It is not difficult to 
verify that the tournament $R_n$ is regular and simple. Moreover, by 
reversing the arc $\left(0,\frac{n-1}{2}\right)$, we obtain a non-simple 
tournament. Hence, the simplicity index of $R_n$ is $1$. If $n$ is prime 
and $n\equiv 3 \pmod 4$, the Paley tournament $P_n$ is also regular but 
its simplicity index is $\frac{n-1}{2}$.
    
	Let $T$ be an $n$-tournament with vertex set $\{v_1, \ldots , v_n\}$. 
The sequences $L_1=(\delta^{+}_{T}(v_i))_{1\leq i\leq n}$
and $L_2=(\delta^{+}_{T}(v_i,v_j))_{1\leq i<j\leq n}$ are frequently used in our study
of the simplicity index. It is natural to ask whether the simplicity index
of $T$ can be expressed in terms of $L_1$ and $L_2$.

\bibliographystyle{plain}
\bibliography{bibpaper}

\end{document}